\documentclass{amsart}

\usepackage{a4wide,amssymb}

\let\d\partial

\def\1{^{-1}}

\def\CP{{\mathbf C\mathbf P}}
\def\MU{{\mathbf M\mathbf U}}
\def\L{{\mathbf L}}
\def\sn{{\operatorname{sn}}}

\newtheorem{theorem}{Theorem}[section]
\newtheorem{lemma}[theorem]{Lemma}

\theoremstyle{definition}

\newtheorem{corollary}[theorem]{Corollary}

\theoremstyle{remark}

\numberwithin{equation}{section}

\begin{document}

\title[]{On addition theorems related to elliptic integrals}

%\date{ Debrecen, May 11, 2017 }

\begin{abstract} We give the formul\ae for the components of the  Buchstaber formal group law and its exponent over $\mathbb{Q}[p_1,p_2,p_3,p_4]$. 
This leads to an addition theorem for the general elliptic integral $\int_0^x dt/R(t)$, for $R(t)=\sqrt{1+p_1t+p_2t^2+p_3t^3+p_4t^4}$. The motivation comes from Euler's addition theorem for elliptic integrals of the first kind.
\end{abstract}
\author{Malkhaz Bakuradze}
\address{Faculty of exact and natural sciences, Iv. Javakhishvili Tbilisi State University, Georgia }
\email{malkhaz.bakuradze@tsu.ge}
\thanks{The first author was supported by CNRS PICS N 7736, and by 
Shota Rustaveli NSF grant 217-614}
\author{Vladimir Vershinin}
\address{University Motpellier 2}
\email{vershini@uni-montp2.fr}
\thanks{Второй автор был поддержан грантом CNRS PICS N 7736}

\maketitle

\section{Introduction and statements }

The Jacobi elliptic sine is the elliptic version of the circular sine. 
In traditional notation, $\sn(u)$ is the inversion of the elliptic integral of the first kind

\begin{align*}
\sn(u)&=x, \,\,\, u=\int_0^x\frac{dt}{\sqrt{(1-t^2)(1-k^2t^2)}}, 
\end{align*}
where $k$ is some parameter called modulus.

Addition theorems offer a means of determining the value of the function for the sum of two quantities as arguments, when the values of the function for each argument is known.

The function $f(u)=\sn(u)$ satisfies Euler's addition theorem  

\begin{equation}
\label{euler}
f(u+v)=\frac{f(u)R(f(v))+f(v)R(f(v))}{1-k^2f(u)^2f(v)^2},\,\,\,\text { where   } R(u)=\sqrt{(1-u^2)(1-k^2u^2)}.
\end{equation}

\medskip

Let us rewrite \eqref{euler} in the form of the addition law due to Cayley \cite{W}, (see also \cite{BU-UST}, (1.4)).

\begin{align}
	\label{addition-ell}
	f(u+v)=&\frac{f^2(u)-f^2(v)}{f(u)f'(v)-f(v)f'(u)},\,\,\, f(0)=0,\,\,\,f'(0)=1.
\end{align}

\bigskip

The further solutions of \eqref{addition-ell}
are as follows. 

\medskip

% The hyperbolic modular sine $\mathcal{SN}(u)$  (see \cite{Ernst})
%is defined by  
%$$
%\mathcal{SN}(u)=x,\,\,\, u=\int_0^x\frac{dt}{\sqrt{(1+t^2)(1+k^2t^2)}}.
%$$ 

%Let $\mathbb{Q}[\delta, \epsilon]$ be a polynomial algebra on two variables %over rational numbers $\mathbb{Q}$.

The exponent series of the Ochanine elliptic genus $f=\exp(\phi_{Oc})$ over  $\mathbb{Q}[\delta, \epsilon]$ is defined as 
the inversion of the elliptic integral 

\begin{equation}
\int_0^x\frac{dt}{\sqrt{1+\delta t^2+\epsilon t^4}}.
\end{equation}
%This is the reason that the Ochanine genus is called the elliptic genus. 

\medskip

\medskip

Let $f_0(u)=1/\Phi(u)$, where  $\Phi(u)=\Phi(u,\alpha)$ is the simplest Baker-Akhiezer function  \cite{KR}. One has the addition formula \cite{BUCH1}

\begin{equation}
f_0(u+v)=
\frac{f_0^2(u)\big(\frac{f_0(v)}{-f_0(-v)}\big)-f_0^2(v)\big(\frac{f_0(u)}{-f_0(-u)}\big)}{f_0(u)f_0'(v)-f_0(v)f_0'(u)}.
\end{equation}

\bigskip

For other interesting examples see \cite{BUCH1}, \cite{BUCH-BU}, \cite{BP}.

%From now on all the series under consideration are over some %$\mathbb{Q}$-algebra.

\medskip

Consider the universal Buchstaber formal group law, the universal example of the formal group law of the form   
\begin{equation}
\label{F}
F(x,y)=\frac{x^2A(y)-y^2A(x)}{xB(y)-yB(x)},
\end{equation}
also specializing to the Euler formal group law  \eqref{euler} for $A(x)=1$ and $B(x)=R(x)$.

It is known from \cite{BUCH1} that the exponent series of \eqref{F}  gives the universal example of $f(u)$ such that $f(0) = 0$, $f'(0) = 1$ and such that $f$ has an addition theorem of the form
\begin{equation}
\label{BUCH-BU}
f(u + v) =\frac{f(u)^2\xi_1(v)-
f(v)^2\xi_1(u)}{f(u)\xi_2(v)-f(v)\xi_2(u)},
\end{equation}
for some series $\xi_1(u)=A(f(u))$ and $\xi_2(u)=B(f(u))=f'(u)$ such that $\xi_1(0) = \xi_2(0) = 1$.
  
%(see proof of Corollary \ref{main-corollary}).

\medskip

By \cite{hoehndiss}  $f$ is defined over the
polynomial ring $\mathbb{Q}[p_1,p_2,p_3,p_4]$ on four variables of degrees 2, 4, 6 and 8.

\medskip

Our main result is Theorem \ref{main}. It provides new explicit formul\ae \eqref{mathcal-B} and \eqref{mathcal-A} for the series $B$ and $A$ in \eqref{F}  after tensoring the coefficient ring of $F$ with rationals, i.e.,  over $\Lambda \otimes \mathbb{Q}=\mathbb{Q}[p_1,p_2,p_3,p_4]$. It also provides the differential equation \eqref{dif} with the general solution $\mathcal{B}$.

%(see the proof of Theorem \ref{theorem-new}.) 

\bigskip

Let us define the following formal power series over the polynomial ring $\mathbb{Q}[p_1,p_2,p_3,p_4]$, where $p_1, \cdots , p_4$ are variables of  degrees 2, 4, 6 and 8 respectively.

\begin{align}
\label{R}
&\mathcal{R}(x):=\sqrt{1+p_1x+p_2x^2+p_3x^3+p_4x^4};\\
\label{mathcal-B}
&\mathcal{B}(x):=
\frac{\mathcal{R}(\mu(x))}{\mu'(x)}
\text{ \,\,\,\,i.e., }\mathcal{R}(x)=
\frac{\mathcal{B}(\nu(x))}{\nu'(x)}  ,\\
\label{nu}
\text { where \,\,\,\,}
&\mu(x):=\frac{x}{\mathcal{B}(x)}, \,\,\,\, \nu(x):=\mu^{-1}(x);
\end{align}

Thus $\mathcal{B}(x)$,  $\mathcal{B}(0)=1$ is a general solution of the differential equation

\begin{equation}
\label{dif}
\mathcal{B}(x)^2(\mathcal{B}(x)-x\mathcal{B}'(x))^2=
\mathcal{B}(x)^4+p_{1} x\mathcal{B}(x)^3+p_2 x^2\mathcal{B}(x)^2+p_3 x^3\mathcal{B}(x)+p_4 x^4.
\end{equation}

\begin{align}
\label{mathcal-A}
&\mathcal{A}(x):=
\mathcal{B}^2(x)-\frac{1}{2}x\mathcal{B}(x)\mathcal{B}'(x)+\frac{1}{4}p_1x\mathcal{B}(x)
-\big( \frac{1}{16}p_1^2-\frac{1}{4}p_2 \big) x^2.
\end{align}

\bigskip

Note that all series above are defined in terms of $\mathcal{R}(x)$.

%The following addition theorem which specializes to Euler's addition theorem %for elliptic integrals of the first kind \eqref{euler}.    

\bigskip

%(and hence for the universal Buchstaber or equivalently Krichever-H\"ohn %formal group laws, see Theorem \ref{theorem-new}.)

%The proof of Theorem \ref{main} is based on the observation that the series %$\nu=(\frac{x}{\mathcal{B}(x)})^{-1}$ is the strict isomorphism from the %formal group law with exponent \eqref{expG} to the formal group law %corresponding to \eqref{expB}. This fact leads to the following  

\bigskip

\begin{theorem}
\label{main} 
In notations \eqref{F} - \eqref{nu} and \eqref{mathcal-A} one has  over the polynomial ring $\mathbb{Q}[p_1,p_2,p_3,p_4]$

\begin{align*}
A(x)=\mathcal{A}(x),\,\,\,B(x)=\mathcal{B}(x),
\text{ and } \log_F(x)=\int_{0}^{x}\frac{dt}{\mathcal{B}(t)}.
\end{align*}

\end{theorem}

\bigskip

Theorem \ref{main} leads to the following

\begin{corollary} 
	\label{main-corollary} Let $f(u)$ be a formal power series over a $\mathbb{Q}$-algebra such that $f(0) = 0$, $f'(0) = 1$. Then $f$ has an addition theorem of the form \eqref{BUCH-BU} if and only if
 it satisfies one of the following properties

 i) $f$ is the inversion of 
	$u=\int_{0}^{x}\frac{dt}{\mathcal{B}(t)}$,
	or equivalently $f$ is a solution of $f'=\mathcal{B}(f)$.

\medskip

ii) The series
$$
h(x):=\frac{f'(x)}{f(x)}
$$
satisfies the differential equation
\begin{equation*}
\label{Hoehn}
(h'(x))^2=S(h(x)),
\end{equation*}
where $S$ is the generic monic polynomial of degree 4:
$$
S(t)=p_4+p_3t+p_2t^2+p_1t^3+t^4.
$$

\end{corollary}

\medskip

Taking into account the remarks after \eqref{BUCH-BU} one has the following	

\begin{corollary}
	\label{xi}
	The following series $\xi_1$ and $\xi_2$ fit into 
\eqref{BUCH-BU}	
	\begin{align}
		\label{xi_1}
		&\xi_1(u)=\mathcal{A}(f(u))=f'^2(u)-
		\frac{1}{2}(f''(u)+f''(0)f'(u))f(u)
		+\frac{1}{2}\big(f''^2(0)-f'''(0)\big)f^2(u);\\
		\label{xi_2}
		& \xi_2(u)=\mathcal{B}(f(u))=f'(u).
	\end{align}

\end{corollary}

\medskip

In \cite{BP} the proof of Theorem E.5.4 derives the  formul\ae which agree  with  \eqref{xi_1}, \eqref{xi_2}. In particular, the coefficient $b_1$ of $x$ in $B(x)$ does not affect on \eqref{F}, and therefore can be chosen arbitrarily. Similarly, the series $\xi_1(u)=A(f(u))$ and $\xi_2(u)=B(f(u))$ in \eqref{BUCH-BU} are defined up to summands $k_1f(u)^2$ and $k_2f(u)$, for any constants $k_1$ and $k_2$ respectively. We derive \eqref{mathcal-A} and \eqref{xi_1} using \eqref{A} of Lemma \ref{Main1}.

\medskip

The following is the addition theorem for general elliptic integral and its inversion $SN(u)$, which specializes to elliptic sine. In particular it is defined by    
\begin{equation}
\label{expG}
SN(u)=x,\,\,\,  u=\int_{0}^{x}\frac{dt}{\mathcal{R}(t)},
\end{equation}
where $\mathcal{R}$ is as in \eqref{R}. The function $SN$ is the exponent series of a genus $\psi $ introduced in \cite{SCH}.

\begin{theorem}
	\label{theorem G(x,y)}
	One has the addition formula

	\begin{equation*}	 
		\int_{0}^{x}\frac{dt}{\mathcal{R}(t)}+\int_{0}^{y}\frac{dt}{\mathcal{R}(t)}
		=\int_{0}^{G(x,y)}\frac{dt}{\mathcal{R}(t)},
	\end{equation*}
	where 
	
	\begin{equation*}
		G(x,y)=\mu \big(P_1+\sigma  P_1+
		\frac{1}{2}\nu(x)\nu(y)\frac{P_2-\sigma P_2}{P_3-\sigma P_3}\big),
	\end{equation*}
	here
	$\sigma$ is the transposition of $x$ and $y$
	and the series $P_i=P_i(x,y)$ are determined by
	\begin{align*}
		&P_1=\nu(x)\mathcal{R}(y)\nu'(y),\\
		&P_2=-\nu(x)\nu'(y)\big(\mathcal{R}(y)
		\big(\mathcal{R}'(y)-\mathcal{R}'(0)\big)
		-\nu(x)\mathcal{R}^2(y)\nu''(y)\big),\\
		&P_3=\nu(x)\mathcal{R}(y)\nu'(y).
	\end{align*}
	
\end{theorem}

\bigskip

\begin{corollary}
	\label{corollary SN}
	Let $SN(x)$ be the inversion of $\int_{0}^{x}\frac{dt}{\mathcal{R}(t)}$ , then  
	\begin{equation*}
		SN(x+y)=G(SN(x),SN(y)).
	\end{equation*}	
\end{corollary}

\bigskip

To prove Theorem \ref{theorem G(x,y)} we first reduce the problem of
 explicit addition theorem to Krichever-H\"ohn genus \cite{hoehndiss}, which
  is defined over $\mathbb{Q}[p_1,p_2,p_3,p_4]$ by \eqref{Hoehn}. For this
   we use the explicit strict isomorphism of Lemma \ref{iso}. From \cite{B1} we know that the universal Buchstaber formal group law $F_B$ \eqref{F} can
    be alternatively defined by the Nadiradze genus $\phi_{N}$ \eqref{phi_N}.  In Theorem \ref{theorem-new} we prove that the Krichever-H\"ohn genus is identical to the Nadiradze genus $\phi_{N}$ after rationalization of its values ring $\Lambda$, i.e., over $\Lambda\otimes \mathbb{Q}=\mathbb{Q}[p_1,p_2,p_3,p_4]$. This reduces the task to the universal Buchstaber formal group law $F_B$ and we apply our explicit
     formul\ae for the components of $F_B$  obtained in Theorem \ref{main}.

\bigskip

\bigskip

 \section{preliminaries}

 It is convenient to give proofs of our results in terms of formal group laws. We give here the necessary definitions and facts.  We refer the reader to \cite{BUCH2} as the detailed survey in the subject.
 
 \medskip
 
 A formal group law  over a commutative ring with unit $R$ is a formal power series $F(x,y)$ in $R[[x,y]]$ satisfying
 
 \medskip
 
 (i) $F(x,0)=F(0,x)=x$,
 
 (ii) $F(x,y)=F(y,x)$,
 
 (iii) $F(x,F(y,z))=F(F(x,y),z)$.
 
 \medskip
 
 Let $F$ and $G$ be formal group laws. A homomorphism from $F$ to $G$ is a power series $h(x)\in R[[x]]$ with constant term $0$ such that
 $$h(F(x,y))=G(h(x),h(y)).$$
 
 It is an isomorphism if $h'(0)$ (the coefficient at $x$) is a unit in $R$, and a strict isomorphism if the coefficient at $x$ is 1.
 
 If $F$ is a formal group law over a commutative $\mathbb{Q}$-algebra $R$, then it is strictly isomorphic to the additive formal group law $x+y$. In other words, there is a strict isomorphism $l(x)$ from $F$ to the additive formal group law. The series $l(x)$ is called the logarithm of $F$, so we have $F(x,y)=l^{-1}(l(x)+l(y))$.  The inverse to logarithm is called the exponential of $F$.

 The logarithm $\log_F(x):=l(x)\in R \otimes \mathbb{Q}[[x]]$ of a formal group law $F$ is given by
 $$
 \log(x)=\int_{0}^{x}\frac{dt}{\omega(t)},\,\,\,\,\,\omega(x)=\frac{\partial F(x,y)}{\partial y}(x,0).
 $$
 
 We will often use the following consequence of definitions above: If  $h(x)$ is a strict isomorphism from the formal group law $F$ to $G$, then
 $$ 
 \log_F=\log_{G}(h),\,\,\,\, \text{ i.e., }  \log_{G}=\log_F(h^{-1}).
 $$

 There is a ring $\L$, called the universal Lazard ring, and a universal formal group law $F(x,y)=\sum a_{ij}x^iy^j$ defined over $\L$. This means that for any formal group law $F'$ over any commutative ring with unit $R$ there is a unique ring homomorphism $r:\L\rightarrow R$ such that $F'(x,y)=\sum r(a_{ij})x^iy^j$.

 The formal group law of geometric cobordism was introduced in \cite{NOV}. Following Quillen we will identify it with the universal Lazard formal group law as it is proved in \cite{Q} that the coefficient ring of complex cobordism $\MU^*=\mathbb{Z}[x_1,x_2,...]$, $|x_i|=2i$ is naturally isomorphic as a graded ring to the universal Lazard ring.

 The coefficients of the formal group law of geometric cobordisms $F_U$ and its logarithm may be described geometrically by the following results.
 
 \medskip
 
 \noindent{\bf Theorem} (Buchstaber \cite{BU}).
 $$
 F_U(u,v)=\sum_{i,j\geq 0}[H_{ij}]u^iv^j(\sum_{r\geq 0}[\CP_r]u^r)(\sum_{s\geq 0}[\CP_s]v^s),$$
 where $\CP_r$ are the complex projective spaces, $H_{ij} (0\leq i \leq j)$ are Milnor hypersurfaces and $H_{ji}=H_{ij}$.
 
 \medskip
 
 \noindent{\bf Theorem} (Mishchenko, see \cite{NOV}). The logarithm of the formal group law of geometric cobordisms is given by the series
 $$\log(u)=u+\sum_{ k\geq 1}\frac{[\CP_k]}{k+1}u^{k+1} \in \MU^*\otimes \mathbb{Q}[[u]].$$

\bigskip

 The addition formula \eqref{euler} corresponds to Euler's addition formula for the elliptic integrals of the first kind

 \begin{equation*}
 %\label{Euler-1}
 \int_0^x\frac{dt}{\sqrt{(1-t^2)(1-k^2t^2)}}+\int_0^y\frac{dt}{\sqrt{(1-t^2)(1-k^2t^2)}}=\int_0^{T(x,y)}\frac{dt}{\sqrt{(1-t^2)(1-k^2t^2)}},
 \end{equation*}
 where
 \begin{equation*}
 %\label{FGL-Euler-1}
 T(x.y)=\frac{x\sqrt{(1-y^2)(1-k^2y^2)}+
 	y\sqrt{(1-x^2)(1-k^2x^2)}}{1-k^2x^2y^2}
 \end{equation*}	
is the Euler formal group law. In terms of the logarithm of the formal group law this means

 \bigskip

 The following two ideas come naturally in mind.
 
 \medskip

 First, one can replace $\sqrt{(1-x^2)(1-k^2x^2)}$ by more general expression \eqref{R}
 and consider the corresponding formal group law $G(x,y)$ with logarithm
 \begin{equation}
 \label{Mog}
 \log_{G}=\int_{0}^{x}\frac{dt}{\mathcal{R}(t)},
 \end{equation}
 specializing to the elliptic integral of the first kind. The corresponding formal group law was recently studied in \cite{SCH}.

 \bigskip

 Second, one can consider the universal Buchstaber formal group law, the universal example of the formal group law of the form  \eqref{F}. 
 
 \medskip

Motivated by string theory, the general complex elliptic genus (also called Krichever-H\"ohn genus) has been defined: in \cite{KR} Krichever wrote down its characteristic power series $Q(x)$ using Baker-Akhiezer function. 
In \cite{hoehndiss}, H\"ohn defined four variable elliptic genus $\phi_{KR}$ determined by the following property: if one denotes by $f=f_{KH}$ the exponent of the corresponding formal group $F_{KH}$, then the series
 $$
 h(x):=\frac{f'(x)}{f(x)}
 $$
 satisfies the differential equation
 \begin{equation}
 \label{Hoehn}
 (h'(x))^2=S(h(x)),
 \end{equation}
 where $S$ is the generic monic polynomial of degree 4,
 $$
 S(t)=p_4+p_3t+p_2t^2+p_1t^3+t^4.
 $$
 
 He also showed that $\phi_{KR}$ takes values in the ring $\mathbb{Q}[p_1,p_2,p_3,p_4]$ and that it agrees with Krichever’s deﬁnition.

 \bigskip

 We use the observation that these two formal groups $G$ and $F=F_{KH}$ are connected by an explicit isomorphism. In particular

 \begin{lemma}
 	\label{iso}
 	The series $\frac{x}{\mathcal{B}(x)}$
is the strict isomorphism from  $F_{KH}$ to  $G$.
 \end{lemma}
 For the proof see \cite{B} or \cite{B2}. Note that on page 13 in \cite{B} the value $\psi(\mathbb{C}P)$ should read as $$\frac{35}{128}q_1^4-\frac{15}{16}q_1^2q_2+\frac{3}{8}q_2^2+\frac{3}{4}q_1q_3-\frac{1}{2}q_4.$$

 \bigskip

 \bigskip

We shall need some results of \cite{B}, \cite{B1} and we collect them together in the following.

Recall here how the universal Nadiradze formal group $F_N$ is constructed \cite{B1}. Let $F_U$ be the universal formal group law over $\MU^*$. Define the series 
$$\sum A_{ij}x^iy^j=F_U(x,y)(x\omega(y)-y\omega(x)).$$
Now kill all $A_{i,j}$, with $i,j\geq 3$. Then the Nadiradze formal group law is classified by the quotient map 

\begin{equation}
\label{phi_N}
\phi=\phi_N:\MU^* \to \MU^*/(A_{i,j}, i,j\geq 3 ):=\Lambda.
\end{equation} 

In other words $F_N$ is the universal formal group law whose invariant differential $\omega(x)\in \Lambda[[x]]$ is 
\begin{equation*}
\frac{\partial F(x,y)}{\partial y}(x,0)=1/\phi(\CP(x))=1/(1+\sum \phi[\CP_i]x^i).
\end{equation*}

Let $F_B$ be the universal Buchstaber formal group law, the universal example among the formal group laws of the form \eqref{F}.

\begin{lemma}  
\label{Main1} \cite{B1} i) The universal Nadiradze formal group law $F_N$ over $\Lambda$ is identical to the universal Buchstaber formal group law $F_B$, i.e., $F_B$ is the formal groups of the form  \eqref{F}, where
	$$
	 B(x)=1/\phi(\CP(x)).
	$$

	ii) The series $A(x)\in \Lambda[[x]]$ can be defined by
	\begin{equation}
	\label{A}
	A(x)=-x^2B(x)\beta(x)-b_1xB(x)+B^2(x)-b_2x^2,
	\end{equation}
	
	where
	$$
	\beta(x)=\frac{B'(x)-b_1}{2x}\in \Lambda[[x]]
	$$
	\medskip
	and $b_1, b_2$ are the coefficients of $B(x)=1+\sum_{i\geq 1}b_ix^i$.

	\end{lemma}
	
	Note that for the formula \eqref{A} we have to use the series $\sum A_{ij}x^iy^j$ in the proof of Proposition 2 in \cite{B1} and rewrite $A(x)=\sum A_{i2}x^i$ and $B(x)=\omega(x)$.

\bigskip

One has also the following 

\begin{theorem}
\label{theorem-new}
Let  $\Lambda$ be the ring of coefficients of the universal Nadiradze formal group $F_N$. Then $F_N$ over $\Lambda \otimes \mathbb{Q}=\mathbb{Q}[p_1,p_2,p_3,p_4]$  is identical to the Krichever-H\"ohn formal group law $F_{KH}$.
\end{theorem}

\bigskip

\section{Proofs}

\subsection*{Proof of Theorem \ref{theorem-new}} By \cite{hoehndiss} the coefficient ring $\Lambda$ of the universal Buchstaber formal group law tensored by rationals is 
$$\Lambda \otimes \mathbb{Q}=\mathbb{Q}[p_1,p_2,p_3,p_4],$$
therefore it suffices to consider the formal group laws over $\mathbb{Q}[p_1,p_2,p_3,p_4]$  

\medskip

 Let $l=log_G$ be the logarithm of the formal group law $G$ in Lemma \ref{iso} and let $j(x)=1/l'$.
 Lemma \ref{iso} says that for a formal group law $F$, the condition of H\"ohn \eqref{Hoehn} is satisfied if and only if
 $j(x)^2$ is a polynomial of degree 4 with constant term 1.

\medskip

Let  $\phi$ be the classifying map of a Buchstaber formal group law and let
%$\CP(x)=1+\sum_{i\geq 1}\CP_ix^i$, where $\CP_i$ is the bordism class %represented by complex projective space of dimension $2i$. Let
		$$
		C(x)=\sum C_ix^{i}=\phi(\CP(x))=1+\sum_{i\geq 1}\phi[\CP_i]x^i.
		$$

 By Lemma \ref{Main1} i) necessarily we have 
	$$
	B(x)=1/C(x),
	$$
	so that we obtain
	$$
	A(x,y):=F(x,y)(x/C(y)-y/C(x))=x^2A(y)-y^2A(x).
	$$
	It follows that
	$$
	\frac{\d^3 A}{\d y^3}(x,0)=cx^2
	$$
	for some constant $c$.
	
	Taking into account that $F(x,y)=f(g(x)+g(y))$, that $g(0)=0$,  $g'(0)=1$, that $f(g(x))=x$, and that $f'(g(x))=1/g'(x)=1/C(x)$ one obtains
	\begin{multline*}
		\frac{\d^3A}{\d y^3}(x,0)
		=
		-\frac{xC''(x)}{C(x)^4}
		+\frac{3xC'(x)^2}{C(x)^5}
		+\frac{3C'(x)}{C(x)^{4}}
		-\frac{3C_{{1}}}{C(x)^{2}}
		+\frac{(3{C_{{1}}}^{2}-4C_{{2}})x}{{C}(x)}\\
		+(-6C_1^3+12C_1C_2-4C_3)x^2
	\end{multline*}
	so that the series $C(x)$ satisfies the differential equation
	$$
	xC(x)C''(x)
	-3xC'(x)^2
	-3C(x)C'(x)
	+3C_1C(x)^3
	-(3{C_{{1}}}^{2}-4C_{{2}})xC(x)^4
	+kx^2C(x)^5=0
	$$
	for some constant $k$.
	
	Let us now substitute in this equation $\nu(x)$ in the place of $x$, where  $\nu(x)$ is inversion of $xC(x)$, i.~e. $\nu(x)C(\nu(x))=x$, so that
	\begin{align*}
		C(\nu(x))&=\frac x{\nu(x)},\\
		C'(\nu(x))&=\nu'(x)\1\frac{\nu(x)-x\nu'(x)}{\nu(x)^2},\\
		\nu(x)C(\nu(x))C''(\nu(x))&=-\nu'(x)^{-3}\frac x{\nu(x)}\nu''(x)
		+2\nu'(x)\1\frac x{\nu(x)}\frac{x\nu'(x)-\nu(x)}{\nu(x)^2}.
	\end{align*}
	We then obtain that $\nu(x)$ satisfies the differential equation
	$$
	2x\nu(x)^2\nu''(x)
	+(kx^3+(6C_1^2-8C_2)x^2-6C_1x-1)x^2\nu'(x)^3
	-2x\nu(x)\nu'(x)^2
	+6\nu(x)^2\nu'(x)=0.
	$$
	Next let us consider  
	$$j(x)=\frac{1}{l'(x)},\,\,\,l=\log_G$$
	mentioned above.
	
%Note $j(x)=\frac{\nu(x)}{x \nu'(x)}$: 

By the strict isomorphism we have $l(x)=g(\nu(x))$. Hence
\begin{align*}
&f(l(x))=f(g(\nu(x)))=\nu(x),\\
&f'(l(x))l'(x)=\nu'(x),\\
\end{align*}
and taking into account
$$\frac{f(l(x))}{f'(l(x))}=1/h((1/h)^{-1}(x))=x$$
%$$\frac{1}{(l'(x)}=\frac{f'(l(x)}{nu'(x)}=\frac{\nu(x)}{x\nu'(x)}$$
we have
\begin{align*}
		j(x)&=\frac{\nu(x)}{x\nu'(x)},\\
		\nu'(x)&=\frac{\nu(x)}{xj(x)},\\
		\nu''(x)&=\frac{\nu(x)(1-j(x)-xj'(x))}{x^2j(x)^2}.
	\end{align*}
	Thus for the series $j$ we obtain the differential equation
	$$
	2xj(x)j'(x)=4j(x)^2+kx^3+(6C_1^2-8C_2)x^2-6C_1x-1.
	$$
	Then further substituting $j(x)=\sqrt{p(x)}$ we obtain
	$$
	xp'(x)=4p(x)+kx^3+(6C_1^2-8C_2)x^2-6C_1x-1
	$$
	and it is easy to see that the general solution of this equation is a fourth degree polynomial.

Conversely,  by \cite{BUCH1} $F_{KR}$ is identical to $F_{B}$ over $\mathbb{Q}[p_1,p_2,p_3,p_4]$ and by Lemma \ref{Main1} i) $F_B$ is identical to $F_N$.

\qed

One can prove  that $F_{KH}$ is a Buchstaber formal group law by Lemma \ref{iso} again. For $\l=\log_G$  we have $\l'=\frac{1}{\mathcal{R}(x)}$ therefore we have for  $\mu(x)=\frac{x}{\mathcal{B}(x)}$, $\nu=\mu^{-1}$ and $g:=\log_{F_{KH}}=l(\mu)$ 
$$g'=(l(\mu))'=\frac{1}{\mathcal{R}(\mu))}\mu':=\frac{1}{\mathcal{B}}.$$
Therefore differentiating the inverse function we get
$$f'=\mathcal{B}(f).$$ 
By Corollary \ref{main-corollary} i) the corresponding formal group law is the Buchstaber formal group law over $\Lambda \otimes \mathbb{Q}=\mathbb{Q}[p_1,p_2,p_3,p_4]$ classified by the obvious inclusion $\Lambda \to \Lambda \otimes \mathbb{Q}$. 
 
 \bigskip

\subsection*{Proof of Theorem \ref{main}}

Let $\l_G$ and $\l_{F_{KH}}$ be the logarithms of the formal group laws $G$ and $F_{KH}$ respectively in Lemma \ref{iso}.
Then because of the strict isomorphism of Lemma \ref{iso} $\mu(x)=\frac{x}{\mathcal{B}(x)}$ one has $\l_{F_{KR}}(x)=\l_{G}(\nu(x))$. By definition  \eqref{Mog} of $\l_{G}$  we have \eqref{mathcal-B}  and 
$\l_{F_{KR}}(x)=\int_{0}^{x}\frac{dt}{\mathcal{B}(t)}$. But the latter coincides with the logarithm of the universal Buchstaber formal group law by Theorem \ref{theorem-new}. 

Now let us prove \eqref{mathcal-A} for the series $A$ and $B$ over  $\mathbb{Q}[p_1,p_2,p_3,p_4]$. By the note after Lemma \ref{Main1} one has  

\begin{align*}
	&\mathcal{A}(x):=
	-x^2\mathcal{B}(x)\beta(x)-b_1x\mathcal{B}(x)+\mathcal{B}^2(x)-b_2x^2,\\
	&\beta(x)=\frac{\mathcal{B}'(x)-\mathcal{B}'(0)}{2x}=
	\frac{\mathcal{B}'(x)-b_1}{2x},\\
\end{align*}

or  

%\begin{align*}
%	&\mathcal{A}(x):=
%-x^2\mathcal{B}(x)\frac{\mathcal{B}'(x)-b_1}{2x}-b_1x\mathcal{B}(x)+
%\mathcal{B}^2(x)-b_2x^2,\\
%\end{align*}

\begin{equation*}
	\label{mathcal-A-2}
	\mathcal{A}(x):=
	\mathcal{B}^2(x)-\frac{1}{2}x\mathcal{B}(x)\mathcal{B}'(x)-
	\frac{1}{2}b_1x\mathcal{B}(x)-b_2 x^2.
\end{equation*}

Therefore it suffices to see

\begin{equation}
\label{b1,b2}
b_1=-\frac{1}{2}p_1;\,\,\,\,b_2=\frac{1}{16}p_1^2-\frac{1}{4}p_2.
\end{equation}

By Lemma \ref{Main1} i) and Theorem \ref{theorem-new} 

\begin{align*}
	\mathcal{B}(x)=&\phi(1/\CP(x))=
	\frac{1}{(1+\phi[\CP_1]x+
		\phi[\CP_2]x^2+\cdots)}=\\
	&1-\phi[\CP_1]x+\phi([\CP_1]^2-[\CP_2])x^2+\cdots\\
\end{align*}

Then there is a formula of \cite{B2} for calculating the values of $\phi=\phi_{KH}$ on $[\CP_i]$, the generators of $\MU^*\otimes \mathbb{Q}$. In particular 

\begin{align*}
	&\phi([\CP_1]= \frac{1}{2}p_1;\\
	&\phi([\CP_1]^2-[\CP_2])= \frac{1}{4}p^2_1-\frac{3}{16}p_1^2-\frac{1}{4}p_2=\frac{1}{16}p_1^2-\frac{1}{4}p_2.
\end{align*}

This proves \eqref{mathcal-A}.
\qed

\bigskip

\subsection*{Proof of Corollary \ref{main-corollary}}  i) follows directly from Theorem \ref{main}. For ii) apply Theorem \ref{theorem-new}.  Here are some comments

Buchstaber's formal group law $F$ is of the form  \eqref{formula-B}. Hence

\begin{equation*} 
	\label{formula-B} 
	f(u+v)=F(f(u),f(v))=\frac{f(u)^2A(f(v))-f(v)^2A(f(u))}{f(u)B(f(v))-f(v)B(f(u))}.
\end{equation*}

Thus $f$ has an addition theorem of the required form and 

\begin{align*}
	&\xi_1(v)=\mathcal{A}(f(u)),\\
	&\xi_2(v)=\mathcal{B}(f(u)).\\
\end{align*}

Let $f(u)$ be a series over any $\mathbb{Q}$-algebra, if it satisfies the addition theorem \eqref{BUCH-BU} and $log=f^{-1}$, then the corresponding formal group law is as follows

$$
F(u,v)=f(\log u+\log v)=\frac{f(\log u)^2 \xi_1(\log v)-f(\log v)^2 \xi_1(\log u)}
{f(\log u)\xi_2(\log v)- f(\log v)\xi_2(\log u)}=\frac{u^2 \xi_1(\log v)-v^2 \xi_1(\log u)}
{u\xi_2(\log v)-v\xi_2(\log u)}.
$$  

By its form $F$ is a Buchstaber formal group law. 

\bigskip

\medskip

\subsection*{Proof of Corollary \ref{xi}}

\eqref{xi_2} is just a differentiation of the inverse function in Theorem \ref{main}.

\eqref{xi_1}: We have $\xi_1(u)=\mathcal{A}(f(u))$.
\begin{align*}
\text{ Let }f(x)=x+f_1x^2+f_2x^3+\cdots.
\end{align*}

Then by \eqref{xi_2} we have $\mathcal{B}(f(x))=f'(x)$, that is
		 
\begin{align*}
&1+b_1(x+f_1x^2+f_2x^3+\cdots)+b_2(x+f_1x^2+f_2x^3+\cdots)^2+\cdots = \\ &1+2f_1x+3f_2x^2+\cdots \\
\end{align*}

hence

\begin{align*}
	 	&b_1=2f_1=f''(0); \\
        &b_2=3f_2-2f_1^2=1/2f'''(0)-1/2f''^2(0).\\
 \end{align*}

Now taking into account  $\mathcal{B}(f(u))=f'(u)$ and
$\mathcal{B}'(f(u))=f''(u)/f'(u)$ we get by \eqref{mathcal-A}

\begin{align*}
	\xi_1(u)=\mathcal{A}(f(u))=&\\
	&\mathcal{B}(f(u))^2-\frac{1}{2}f(u)\mathcal{B}(f(u))\mathcal{B}'(f(u))
	-\frac{1}{2}b_1f(u)\mathcal{B}(f(u))
	-b_2 f(u)^2= \\
	& f'^2(u)-\frac{1}{2}f(u)f''(u)-\frac{1}{2}b_1f(u)f'(u)
	-b_2f(u)^2=\\
	&f'^2(u)-\frac{1}{2}(f''(u)+f''(0)f'(u))f(u)+\frac{1}{2}\big(f''^2(0)-f'''(0)\big)f^2(u).\\
	\end{align*}

\qed

\bigskip

\subsection{Proof of Theorem \ref{theorem G(x,y)}}

i) By Lemma \ref{iso} one has

\begin{equation*}
	G(x,y)=\nu^{-1} \big[\frac{\nu(x)^2\mathcal{A}(\nu(y))-\nu(y)^2\mathcal{A}(\nu(x))}
	{\nu(x)\mathcal{B}(\nu(y))-\nu(y)\mathcal{B}(\nu(x))} \big]. 
\end{equation*}

\medskip

Then \eqref{mathcal-B} and $\mu=\nu^{-1}$ imply

\begin{align*}
&\mathcal{B}(\nu(x))=\frac{\mathcal{R}(\mu(\nu(x))}{\mu'(\nu(x))}=\mathcal{R}(x)\nu'(x),\\
&\mathcal{B}'(\nu(x))=\frac{\mathcal{B}(\nu(x))'}{\nu'(x)}=
\frac{(\mathcal{R}(x)\nu'(x))'}{\nu'(x)}=
\frac{\mathcal{R}'(x)\nu'(x)+\mathcal{R}(x)\nu''(x)}{\nu'(x)}.
\end{align*}

\bigskip

By \eqref{mathcal-A} we have

\begin{align*}
	\mathcal{A}(\nu(x))=\\
	&\mathcal{R}^2(x)(\nu'(x))^2-\frac{1}{2}\nu(x)\mathcal{R}(x)\nu'(x)\frac{\mathcal{R}'(x)\nu'(x)+\mathcal{R}(x)\nu''(x)}{\nu'(x)}+\\
	&\frac{1}{4}p_1\nu(x)\mathcal{R}(x)\nu'(x)
	-\big(\frac{1}{16}p_1^2-\frac{1}{4}p_2 \big) \nu^2(x)=\\
	&\mathcal{R}^2(x)(\nu'(x))^2-
	\frac{1}{2}\nu(x)\mathcal{R}(x)\mathcal{R}'(x)\nu'(x)
	+\frac{1}{4}p_1\nu(x)\mathcal{R}(x)\nu'(x)\\
	-&\frac{1}{2}\nu(x)\mathcal{R}^2(x)\nu''(x)-\big(\frac{1}{16}p_1^2-\frac{1}{4}p_2 \big) \nu^2(x)=\\
	&\mathcal{R}^2(x)(\nu'(x))^2-
	\frac{1}{2}\nu(x)\mathcal{R}(x)\nu'(x)\big[\mathcal{R}'(x)-\mathcal{R}'(0)\big]\\
	-&\frac{1}{2}\nu(x)\mathcal{R}^2(x)\nu''(x)
	-\big(\frac{1}{16}p_1^2-\frac{1}{4}p_2 \big) \nu^2(x)
\end{align*}
as $\mathcal{R}'(0)=\frac{1}{2}p_1$.

\medskip

Now note that

\begin{align*}
	&\frac{\mathcal{A}(\nu(y))\nu^2(x)-\mathcal{A}(\nu(x))\nu^2(y)}
	{\nu(x)\mathcal{R}(y)\nu'(y)-\nu(y)\mathcal{R}(x)\nu'(x)}=
	\nu(x)\mathcal{R}(y)\nu'(y)+\nu(y)\mathcal{R}(x)\nu'(x)+
	\frac{I+II}{III},
\end{align*}
where

$$
I=
-\frac{1}{2}\nu^2(x)\nu(y)\nu'(y)\mathcal{R}(y)
\big[\mathcal{R}'(y)-\mathcal{R}'(0)\big]
+\frac{1}{2}\nu^2(y)
\nu(x)\nu'(x)\mathcal{R}(x)\big[\mathcal{R}'(x)-\mathcal{R}'(0)\big],
$$

$$
II=
-\frac{1}{2}\nu^2(x)\nu(y)\mathcal{R}^2(y)\nu''(y)
+\frac{1}{2}\nu^2(y)\nu(x)\mathcal{R}^2(x)\nu''(x),	
$$

$$
III=\nu(x)\mathcal{R}(y)\nu'(y)-\nu(y)\mathcal{R}(x)\nu'(x).
$$

\medskip

Thus we get

\begin{equation*}
	\nu(G(x,y))=\nu(x)\mathcal{R}(y)\nu'(y)+\nu(y)\mathcal{R}(x)\nu'(x)+
	\frac{1}{2}\nu(x)\nu(y)\frac{(P-\sigma P)}{\nu(x)\mathcal{R}(y)\nu'(y)-\nu(y)\mathcal{R}(x)\nu'(x)},
\end{equation*}
where
$\sigma \in S_2$ is the transposition of $x$ and $y$
and
$$
P=-\nu(x)\nu'(y)\big(\mathcal{R}(y)
\big(\mathcal{R}'(y)-\mathcal{R}'(0)\big)
-\nu(x)\mathcal{R}^2(y)\nu''(y)\big).
$$

This proves Theorem \ref{theorem G(x,y)}.

Finally, as usual, the formal group law $F=G$ gives the addition formula for its exponent $f=SN$
$$f(u+v)=F(f(u),f(v)),$$
and Corollary \ref{corollary SN} follows.

 \end{document}